\documentstyle[12pt]{article}

\begin{document}
\newcommand{\ol }{\overline}
\newcommand{\ul }{\underline }
\newcommand{\ra }{\rightarrow }
\newcommand{\lra }{\longrightarrow }
\newcommand{\ga }{\gamma }
\newcommand{\st }{\stackrel }
\newcommand{\scr }{\scriptsize }
\title{\Large\textbf{ Generalized Covering Groups and \\
 Direct Limits }}
\author{\textbf{
Behrooz Mashayekhy\footnote{Corresponding author: mashaf@math.um.ac.ir} and Hanieh Mirebrahimi} \\
Department of Mathematics,\\ Center of Excellence in Analysis on Algebraic Stuctures,\\
Ferdowsi University of Mashhad,\\
P. O. Box 1159-91775, Mashhad, Iran.}
\date{ }

\maketitle
\begin{abstract}
M. R. R. Moghaddam (Monatsh. Math. 90 (1980) 37-43.) showed that the
Baer invariant commutes with the direct limit of a directed system
of groups. In this paper, using the generalization of Schur's
formula for the structure of a $\mathcal{V}$-covering group for a
Schur-Baer variety $\mathcal{V}$, we show that the structure of a
$\mathcal{V}$-covering group commutes with the direct limit of a
directed system, in some senses. It has a useful application in
order to extend some known structures of $\mathcal{V}$-covering
groups for several famous products of finitely many  to an arbitrary
family of groups.
\end{abstract}
A.M.S. Classification: 20E10, 20E18, 20J15, 20F18.\\
\textit{Keywords}: $\mathcal{V}$-covering group,
Direct limit, Baer invariant, Variety of groups. \\
\newpage
\begin{center}
\hspace{-0.65cm}\textbf{1. Introduction }\\
\end{center}

Historically, there have been several papers from the beginning of
the twentieth century trying to find some structures for the
well-known notion the covering group and its varietal
generalization, the $\mathcal{V}$-covering group of some famous
groups and products of groups, such as the direct product, the
nilpotent  and the regular product [4, 5, 8, 11, 15, 21, 24].

It is known that any group has at least a covering group [21, 25].
Also the number of covering groups for an arbitrary group has been
studied by I. Schur [22]. Moreover, it is proved that any group has
a $\mathcal{V}$-covering group, where $\mathcal{V}$ is a Schreier
variety [8, 10, 15].

In 1971, J. Wiegold [24] found a structure of a covering group for
the direct product $G=A\times B$ so that the second nilpotent
product of the covering groups $A^*$ and $B^*$ is a covering group
for $G$. In 1972, W. Haebich [4] constructed a covering group for
a finite regular product. Moreover the structure of a covering
group for the verbal wreath product of two groups has been
studied by W. Haebich [5], in 1977.

Naturally, it is of interest to know which class of groups does {\it
not} have a ${\cal V}$-covering group. The first author [9, 10] gave
some examples of groups which do not have any generalized covering
group with respect to the variety of nilpotent groups of class at
most $c\geq 2$, ${\mathcal{N}}_c$. More precisely, he [10] proved
that every nilpotent group of class $n$ with nontrivial
$c$-nilpotent Schur multiplier does not have any
${\mathcal{N}}_c$-covering group for $c>n$.  Thus the results of
Wiegold and Haebich mentioned above cannot be generalized to an
arbitrary variety.
 Moreover the first author [10] has given a
complete answer to the existence of ${\mathcal{N}}_c$-covering group
for finite abelian groups. Also in 2003, he in a joint paper [11]
found a structure for the ${\mathcal{N}}_c$-covering group of a
nilpotent product of a family of cyclic groups.

Now in this paper, we intend to prove that the structure of
$\mathcal{V}$-covering group, commutes with the direct limit of a
directed system of groups in some senses (see Theorem 3.5).
Furthermore, we give an example showing that the hypothesis of being
directed for the system of groups, is an essential condition (see
Example 3.6). Finally, as an application, we extend some of the
previous formulas for the structure of $\mathcal{V}$-covering groups
for
 finite direct and regular
products of groups to infinite ones.
\begin{center}
\textbf{2. Notation and Preliminaries}
\end{center}

We shall assume that the reader is familiar with the notion of the
verbal subgroup $V(G)$, and the marginal subgroup $V^*(G)$ of a
group $G$, associated with a variety $\mathcal{V}$. Whenever
varieties of groups are discussed, we refer to H. Neumann [17] for
notation and
basic results.\\
\textbf{Definition 2.1.} Let $\mathcal{V}$ be a variety of groups
defined by the set of laws $V$, and let $G$ be a group with a
free presentation
\[1\rightarrow R\rightarrow F\rightarrow G\rightarrow 1,\]
where $F$ is a free group. Then the Baer invariant of $G$, denoted
by $\mathcal{V}M$(G), is defined to be
\[R\cap V(F)/[RV^*F],\]
where $V(F)$ is the verbal subgroup of $F$ and $[RV^*F]$ is the
least normal subgroup $N$ of $F$ contained in $R$ such that
$R/N\subseteq$ $V^*(F/N)$. Thus it is the subgroup generated by
the following set:
\[ [RV^*F]=<v(f_1,\ldots ,f_{i-1},f_ir,f_{i+1},\ldots ,f_n)v(f_1,\ldots
,f_i,\ldots ,f_n)^{-1}\mid \] \[ r\in R\ ,\ f_i\in F,\ v\in V\ ,\
1\leq i\leq n\ ,\ n\in {\bf N}>\ .\]

Note that the Baer invariant of a group $G$ is always abelian and
independent of the choice of the presentation of $G$ and it might be
possible to regard ${\cal V}M(-)$ as the first left derived functor
of the functor from all groups to $\cal V$ taking $G$ to $G/V(G)$
[8]. In particular, if $\mathcal{V}$ is the variety of abelian
groups, then the Baer invariant of the group $G$ will be $R\cap
F'/[R,F]$, which is the well-known notion, the Schur multiplier of
$G$, $M(G)$. Also, if $\mathcal{V}=\mathcal{N}$$_c$ is the variety
of nilpotent groups of class at most $c$, then the Baer invariant of
the group $G$ will be $R\cap \gamma_{c+1}(F)/[R,\ _cF]$, where
$\ga_{c+1}(F)$ is the $(c+1)$-st term of the lower central series of
$F$ and $[R,\ _1F]=[R,F]$, $[R,\ _cF]=[[R,\ _{c-1}F],F]$
inductively. It is also called the
$c$-nilpotent multiplier of $G$.\\
\textbf{Definition 2.2.} A variety $\cal V$ is called a {\it
Schur-Baer variety} if for any group $G$ for which the marginal
factor group $G/ V^*(G)$ is finite, then the verbal subgroup
$V(G)$ is also finite and $|V(G)|$ divides a power of $|G/
V^*(G)|$.

 I. Schur [21] proved that the variety of abelian group is a Schur-Baer variety
and R. Baer [2] proved that a variety defined by some outer
commutator words, for instance the variety ${\mathcal{N}}_c$, has
the above property.
 The following theorem gives a useful property of Schur-Baer
 varieties.\\
\textbf{\bf Theorem 2.3} ([8]).
 A variety $\cal V$ is a Schur-Baer variety if and only if for
 every finite group $G$, its Baer invariant ${\cal V}M(G)$ is of
order dividing a power of $|G|$.\\
\textbf{\bf Definition 2.4.} Let $\cal V$ be a variety of groups
and let $G$ be a group. Then, by definition, a $\cal V$-{\it
covering group} of $G$ (a generalized covering group of $G$ with
respect to $\cal V$) is a group $G^*$ with a normal subgroup $A$
such that $G^*/A\cong G$, $A\subseteq V(G^*)\cap V^*(G^*)$, and
$A\cong {\cal V}M(G)$ (see [8]).

 Note that if $\cal V$ is the variety of abelian groups, then a $\cal
V$-covering group of $G$ will be an ordinary covering group
(sometimes it is called a representing  group) of $G$. Also if
${\cal V}={\cal N}_c$, then an ${\cal N}_c$ -covering group of $G$
is a group $G^*$ with a normal subgroup $A$ such that
$$\frac{G^*}{A}\cong G,\ A\subseteq Z_c(G)\cap {\ga}_{c+1}(G),\ {\it and}\
A\cong {\cal N}_cM(G)\ .$$ It is well-known that every group $G$
has a covering group (see [7, 21, 25]). In general, for the
existence of $\mathcal{V}$-covering groups, we have the following
concepts and results. Let $\cal V$ be a variety of groups, then a
group $G$ is called $\cal V$-free if it is a free object in the
category of all groups in the variety $\cal V$. It is known that
if $F$ is a free group, then $F/ V(F)$ is a $\cal V$-free group.
By a well-known theorem of Schreier, every subgroup of a free
group is also
free. Thus it is natural to define the following notion.\\
\textbf{\bf Definition 2.5.} Let $\cal V$ be a variety of
groups.  Then $\cal V$ is called a Schreier variety if and only
if every subgroup of a $\cal V$-free group is also $\cal V$-free.
It has been proved by P. M. Neumann and J. Wiegold [8, 18] that
the only Schreier varieties are the variety of all groups $\cal
G$, the variety of abelian groups, $\cal A$, and the variety of
all abelian group of exponent p, ${\cal A}_p$, where p is a prime.

 Note that the notion of Schreier varieties can be generalized to varieties
in which every subgroup of a $\cal V$-free group is a $\cal
V$-splitting group i.e. a group $G$ in the variety $\cal V$ which
splits every short exact sequence in $\cal V$ of the following
form
$$1 \ra A \ra B \ra G \ra 1\ .$$
Clearly every Schreier variety has the above property. In fact the
only varieties with the above property are $\cal G$, $\cal A$,
${\cal A}_m$ where m is square free (see [8, 18]).

 C. R. Leedham-Green and S. McKay [8], by a homological method, proved that a
sufficient condition for the existence of a $\cal V$-covering
group of $G$ is that $G/V(G)$ should be a $\cal V$-splitting
group. Also in the following theorem we give a similar sufficient
condition for the existence of a ${\cal V}$-covering group of a
group.\\
\textbf{\bf Theorem 2.6} ([8, 10, 15]). Let $\cal V $ be a variety
in which every subgroup of a $\cal V$-free group is $\cal
V$-splitting. Then every group has a $\cal V$-covering group. In
particular, if $\cal V$ is a Schreier variety, then every group has
a  $\cal V$-covering group.

Also the first author [10] showed that if $G$ is a nilpotent group
of class $n$ such that ${\cal N}_c$$M(G)\neq 1$, and $c>n$, then $G$
has no ${\cal N}_c$-covering group. Moreover, this fact has been
extended to the variety of polynilpotent groups in some senses [12].

We next attend to some concepts as well-known categorical objects,
the direct system and the direct limit which are defined [19,20]
for any arbitrary category, if any. Note that in some famous
categories, specially the category of all groups, every direct
system has a direct limit (see [19,20]). In particular, in the
category of all groups, the free product $\prod^*_{i\in I}G_i$ as
a coproduct of $\{G_i\}_{i\in I}$ is a direct limit of this family
which is the trivial direct system and so is not directed (see
[20]).

In our main results we deal with the particular case of direct
systems which are indexed with a directed set, called directed
system. In fact, we use one of the equivalent form of the
definition of direct limit of a directed system in the category of
all groups, which is more simple and useful for our goal, as we
mention in following preliminaries. Note that in all following
points, we refer the reader to
[14, 19, 20] for further details. \\
\textbf{\bf Definition 2.7.} Let $\{G_i\}$ be a direct system of
groups indexed by a partially ordered set I, which is also
directed, that is, for every $i,j\in I$ there exists $k\in I$
such that $i,j\leq k$. For $i\leq j$, let there exists a
homomorphism $\lambda_i^j:G_i\rightarrow G_j$ such that:\\
$(i)$ $\lambda_i^i:G_i\rightarrow G_i$ is the identity map of
$G_i$,
for all $i\in I$;\\
$(ii)$ If $i\leq j\leq k$, then
$\lambda_i^j\lambda_j^k=\lambda_i^k$, as the following commutative
diagram:
\[\begin{array}{ccccccc}
 G_i  \st{\lambda_i^j}{\longrightarrow} G_j  \\
 \vspace{-0.5cm}\\
 \hspace{0cm}  \st{{\lambda}^k_i}{}\searrow  \hspace{0.4cm} \downarrow \st{\lambda^k_j}{}\\
 \hspace{0.8cm}G_k.\\
\end{array}\]
In this case, we call the system $\{G_i; \lambda_i^j,I\}$ a {\it
directed system}. Now we define an equivalence relation on the
disjoint union $\bigcup_{i\in I}G_i$,  by: if $x\in G_i$ and
$y\in G_j$, then
\begin{center}
$x\sim y$ if and only if $x\lambda_i^k=y\lambda_j^k$ for $k\geq
i,j$.
\end{center}
Let $G$ denote the quotient set $\bigcup_{i\in I}G_i/\sim$ and use
$\{x\}$ for the equivalence class of $x$. Also we define a
multiplication on $G$ as follows: if $\{x\}$, $\{y\}$ are elements
of G, we choose $i,j\in I$ such that $x\in G_i$ and $y\in G_j$
then
\[\{x\}\{y\}=\{(x\lambda_i^k)(y\lambda_j^k)\},\ for\ k\geq i,j.\]
Clearly this is a well-defined multiplication, which makes $G$
into a group and it is called the {\it direct limit} of the
directed system $\{G_i;\lambda_i^j,I\}$. It will be denoted by
\[\ \ \ \ \lim_{\longrightarrow} G_i=\bigcup_{i\in I}G_i/\sim=G.\]

We need only the following well-known results of direct limits.\\
\textbf{\bf Lemma 2.8.} Suppose that $\{G_i; \lambda_i^j,I\}$ is a
directed system of groups and
$G=\displaystyle{\lim_{\longrightarrow} G_i}$. Then we have the
following
statements:\\
$(i)$ The group $G$ has the universal property so that for a given
group $H$ and homomorphisms $\tau_i:G_i\rightarrow H$, such that
$\lambda_i^j\tau_j=\tau_i$ for all $i\leq j$, there exists a
unique homomorphism $\tau:G\rightarrow H$ such that all the
diagrams
\[\begin{array}{ccccccc}
 G_i    \\
  \hspace{0.1cm}\st{\lambda_i}{}\downarrow  \hspace{0.4cm} \searrow \st{{\tau}_i}{}   \\
\vspace{-0.5cm}\\
 \hspace{0.5cm}G\st{\tau}{\longrightarrow}H\\
\end{array}\]
commute, that is $\lambda_i\tau=\tau_i$, for all $i\in I$.\\
$(ii)$ Direct limit of exact sequences, indexed by a directed set,
is exact, and so in this case, the direct limit preserves
injections.\\
$(iii)$ Let $G$ be an arbitrary group, then $G$ is the direct
limit of its finitely generated subgroups, under the obvious
directed system arising from the inclusion maps.\\
\textbf{\bf Definition 2.9.} Let $\mathcal{C}$ and $\mathcal{D}$
be two categories and let $T_1:\mathcal{C}\rightarrow
\mathcal{D}$, $T_2:\mathcal{D}\rightarrow \mathcal{C}$ be two
functors such that for any $X\in \mathcal{C}$, $Y\in \mathcal{D}$
there is a natural equivalence
\[Hom_{\mathcal{D}}(T_1 X,Y)\simeq Hom_{\mathcal{C}}(X,T_2 Y).\]
In this case, $T_2$ is called a {\it right adjoint functor} to
$T_1$ and the pair ($T_1,T_2$) is called an {\it adjoint pair}. It
is well-known fact  that every functor which has a right-adjoint,
commutes with direct limits. So we have the
following lemma (see [20]). \\
\textbf{\bf Lemma 2.10.} Let $\{X_i; \lambda_i^j,I\}$ be a direct
system of sets indexed by a partially ordered set I, and let
$\mathcal{C}$ and $\mathcal{G}$ denote the categories of sets and
groups, respectively. If $F:\mathcal{C}\rightarrow \mathcal{G}$
is the free functor which associates with every set the free
group on that set as basis, then $F$ commutes with direct limit,
that is,
\[F(\lim_{\longrightarrow} X_i)=\lim_{\longrightarrow} F(X_i).\]
\textbf{\bf Note 2.11.} As a corollary, suppose that $\{G_i;
\lambda_i^j,I\}$ is a directed system of groups, and the sequence
\[1\rightarrow R_i\rightarrow F_i\rightarrow G_i\rightarrow 1\]
is a free presentation for $G_i$, where $F_i(=F(G_i))$ is the free
group on the underlying set of $G_i$, for all $i\in I$. Now using
lemma 2.8.(ii), the direct limit of a directed set is an exact
functor, and hence kernel-preserving, so by Lemma 2.10, the sequence
\[1\rightarrow \lim_{\longrightarrow} R_i\rightarrow \lim_{\longrightarrow}
 F_i\rightarrow \lim_{\longrightarrow} G_i\rightarrow 1\]
is a free presentation for $\displaystyle{\lim_{\longrightarrow} G_i}$.\\
\textbf{\bf Lemma 2.12.} With the above assumption and
notation, we have the following relations:\\
$(i)$ $\displaystyle{(\lim_{\longrightarrow} R_i)\cap
V(\lim_{\longrightarrow}
F_i)=\lim_{\longrightarrow} (R_i\cap F_i)}$;\\
$(ii)$ $[(\displaystyle{\lim_{\longrightarrow} R_i) V^*
(\lim_{\longrightarrow} F_i)]=\lim_{\longrightarrow}
[R_i V^* F_i]}$.\\
\textbf{\bf Theorem 2.13} ([14]). Let $\{G_i; \lambda_i^j,I\}$ be a
directed system of groups. Then for a given variety $\mathcal{V}$,
the Baer invariant commutes with direct limit, that is,
\[{\mathcal{V}}M(\lim_{\longrightarrow} G_i)=\lim_{\longrightarrow}
{\mathcal{V}}M(G_i).\]\\

\begin{center}
\textbf{3. The Main Result}
\end{center}

In order to deal with $\cal V$-covering groups of a group $G$, it
is useful to know more relationship between the Baer invariant
${\cal V}M(G)$ and the $\cal V$-covering groups of $G$. In this
aspect, to prove our main theorem, first of all we need to point
the following notes which are the generalization of
some parts of an important theorem of Schur [7, Theorem 2.4.6].\\
\textbf{\bf Lemma 3.1.}\textit{ Let $\cal V$ be a variety of
groups and $G$ be a group with a free presentation $1\ra R\ra
F\ra G\ra 1$. If $S$ is a normal subgroup of $F$ such that
\[ \frac {R}{[RV^*F]}\cong\frac {R\cap V(F)}{[RV^*F]}\times \frac
{S}{[RV^*F]},\ \ \ \ (*)\]
then $G^*=F/S$ is a $\cal V$-covering group of $G$.}\\
\textit{Proof.} Setting $A=R/S$, so $G^*/A \cong F/R\cong G$ and
using (*), we have
\[ A=\frac {R}{S}\cong \frac
{R/[RV^*F]}{S/[RV^*F]}\cong \frac {R\cap V(F)}{[RV^*F]}\cong {\cal
V}M(G).\] Since $[RV^*F]\subseteq S$, therefore $A=R/S\subseteq
V^*(F/S)=V^*(G^*)$. Also we have
$$A=\frac {R}{S}=\frac {RS}{S}\ \st {by(*)}{\subseteq }\ \frac {V(F)S}{S}=V(\frac
{F}{S})=V(G^*)\ .$$
Hence $G^*$ is a $\cal V$-covering group of $G$. $\Box$\\
\textbf{\bf Lemma 3.2.}\textit{ Let $\cal V$ be a variety of
groups, let $G$ be a group with a free presentation $1\ra R\ra
F\ra G\ra 1$, and let $G^*$ be a $\cal V$-covering group of
$G$. Then $G^*$ is a homomorphic image of $F/[RV^*F]$.} \\
\textit{Proof.} Let $F$ be free on $X$ and $\pi :F\ra G$ be an
epimorphism such that $R=ker(\pi )$. Since $G^*$ is a $\cal
V$-covering group of $G$, we have the following exact sequence
$$1\ra A\lra G^*\st {\Phi }{\lra} G\ra 1\ ,$$
where $A\subseteq V^*(G^*)\cap V(G^*)$ and $A\cong {\cal V}M(G)$.
Since $\Phi $ is surjective, there exists $l_x$ in $G^*$ such
that $ \Phi (l_x)=\pi (x)$, for all $x\in X$. Therefore we have
$$G^*=<A,l_x|x\in X>=AN,$$
where $N=<l_x|x\in X>$. Now, by a result of N. S. Hekster [6],
$$ A\subseteq V(G^*)=V(AN)=V(N)[AV^*G^*]\subseteq V(N)\subseteq N \ .$$
(Note that since $A\subseteq V^*(G^*)$, we have $[AV^*G^*]=1$.)
 Hence we have
$$G^*=N=<l_x|x\in X>.$$
Now consider the homomorphism $\Psi :F\ra G^*$ defined by $\Psi
(x)=l_x$, $x\in X$. Then $\Psi $ is surjective and $\pi =\Phi
\circ \Psi $. Since $1=\pi (R)=\Phi (\Psi (R))$, we have $\Psi
(R)\subseteq A$, so that
$$ \Psi ([RV^*F])\subseteq [\Psi (R)V^*G^*]\subseteq [AV^*G^*]=1\ .$$
It follows that $\Psi $ induces an epimorphism $\ol {\Psi }:
F/[RV^*F]\lra
G^*$. $\Box$ \\
\textbf{\bf Theorem 3.3.}\textit{ Let $\cal V$ be a
 Schur-Baer variety and $G$ be a
finite group with a free presentation $1\ra R\ra F\ra G\ra 1$. If
$G^*$ is $\cal V$-covering group of $G$, then there exists a
normal subgroup $S$ of $F$ such that
$$\frac {R}{[RV^*F]}\cong\frac {R\cap V(F)}{[RV^*F]}\times \frac {S}{[RV^*F]}\ ,$$
and so $G^*\cong F/S$.}\\
\textit{Proof.} By the proof of Lemma 3.2 and its notation, for
every $a\in A$, there exists $x\in F$ such that $a=\Psi (x)$.
Hence $\pi (x)=\Phi \circ \Psi (x)=\Phi (\Psi (x))=\Phi (a)=1$, so
$x\in R$ and thus $A\subseteq \Psi (R)$. Also in the proof of
Lemma 3.2 we showed that $\Psi (R)\subseteq A$, so $A=\Psi (R)$.
Next we observe that
$$\Psi (R\cap V(F))\subseteq \Psi (R)\cap \Psi (V(F))= A\cap V(G^*) =A\ .$$
To prove the other inclusion, assume that $z=\Psi (x)=\Psi (y)$ for
some $x\in V(F)$ and $y\in R$. Then $x^{-1}y\in Ker\Psi$, so $\pi
(x^{-1}y)=1$ and therefore $x^{-1}y\in R$. It follows that $x\in R$
and $z\in \Psi (R\cap V(F))$. Thus  $$ A=\Psi (R\cap V(F))\ .$$
 Now $\ol {\Psi}$ defines an epimorphism $$ {\ol {\Psi}_1}:\frac
{R\cap V(F)}{[RV^*F]}\lra A\ .$$ Since $\cal V$ is a Schur-Baer
variety, $A\cong {\cal V}M(G)=R\cap V(F)/[RV^*F]$, and $G$ is
finite, so by Theorem 1.3 $A$ is also finite. Thus the above
epimorphism is isomorphism. Now, put $S=Ker(\Psi )\cap R$. Then
clearly $S\unlhd F$ and $S/[RV^*F]$ is the kernel of the
restriction of $\ol {\Psi}$ to $R/[RV^*F]$, i.e.
$$ {\ol {\Psi}}_2:\frac {R}{[RV^*F]}\lra A\ .$$
 Now we can consider the following short exact sequence
$$1\ra \frac {S}{[RV^*F]}\lra \frac{R}{[RV^*F]}\st{{\ol {\Psi}_2}}{\lra } A\ra 1.\ \ \ (**)$$
Since $\ol {\Psi}_1$ is an isomorphism, so we have ${\ol
{\Psi}_2}\circ {\ol {\Psi}_1}^{-1}=1_A$, and hence the above short
exact sequence splits. Therefore we have
$$\frac {R}{[RV^*F]}\cong A\rhd\!\!\!< \frac {S}{[RV^*F]}\cong \frac
{R\cap V(F)}{[RV^*F]}\rhd\!\!\!< \frac {S}{[RV^*F]}\ .$$ But
clearly $R\cap V(F)/[RV^*F]\unlhd R/[RV^*F]$, so the above
semidirect product is actually a direct product. Now, by Lemma
3.1, $F/S$ is a $\cal V$-covering group of $G$.

 Let $\theta :F/S \lra G^*$ be the homomorphism induced by $\Psi $. Since $\Psi$
is surjective and $\Psi (R)=A$, $\theta $ is surjective and
$\theta (R/S)=A$. However $|F/S|=|G^*|$, so $\theta $ becomes an
isomorphism i.e. $G^*\cong F/S$. $\Box$

Note that this generalization of the Schur Theorem, has been posed
and proved by M. R. R. Moghaddam and A. R. Salemkar [16], but it
seems that there are some missing points in their proof,
specially the splitting of the exact sequence (**), and so the
condition of being Schur-Baer for the variety $\cal V$.\\
\textbf{\bf Lemma 3.4.}\textit{ The direct limit with a directed
index set, as we mentioned in Definition 2.5, preserves the finite
direct product, that is, for any two directed systems of groups
$\{A_i;\lambda_i^j,I\}$ and $\{B_i;\mu_i^j,I\}$, we have}
\[\lim_{\longrightarrow}(A_i\times B_i)=\lim_{\longrightarrow}A_i\times
\lim_{\longrightarrow}B_i.\] \textit{Proof.} Firstly, for any
$i\in I$, we have the following split exact sequence with natural
homomorphisms
\[1\rightarrow A_i\rightarrow A_i\times B_i\rightarrow B_i\rightarrow 1.\]
Now using Lemma 2.8(ii), the direct limit preserves exactness and so
we have the following exact sequence which is also split
\[1\rightarrow \lim_{\longrightarrow} A_i\rightarrow
\lim_{\longrightarrow}(A_i\times B_i)\rightarrow
\lim_{\longrightarrow}B_i\rightarrow 1.\] On the other hand, we know
that the direct limit is kernel-preserving and so preserves normal
subgroups. Hence $\displaystyle{\lim_{\longrightarrow}B_i}$ is a
normal subgroup of
$\displaystyle{\lim_{\longrightarrow}(A_i\times B_i)}$ and so the result holds.\ $\Box$\\

Now, in order to state and prove the main result of the paper, we
need to explain the concept of an induced directed system of $\cal
V$-covering groups which we use in the main theorem. Let $\cal V$
be a Schur-Baer variety, and let $\{G_i; \lambda_i^j,I\}$ be a
directed system of finite groups. suppose that $G_i^*$ is a $\cal
V$-covering group for $G_i$, for all $i\in
I$. Now if we consider the sequence    \\
\[1\rightarrow R_i\rightarrow F_i\rightarrow G_i\rightarrow 1\] as
a free presentation, then using Theorem 3.3 there exists a normal
subgroup $S_i$ of $F_i$ in such a way that $G_i^*=F_i/S_i$ and
specially satisfies the following relation:
\[\frac {R_i}{[R_iV^*F_i]}\cong\frac {R_i\cap V(F_i)}{[R_iV^*F_i]}\times
\frac {S_i}{[R_iV^*F_i]}\ .\ \ \ \ (*)\] By these notations, for
any $i\leq j$ in I, there exists an induced homomorphism
$\widehat{\lambda}_i^j$ commutes the following diagram
\[\begin{array}{ccccccc}
1 \longrightarrow  R_i  \longrightarrow  F_i  \longrightarrow  G_i  \longrightarrow  1\\
 \hspace{1.6cm}  \downarrow \st{{\widehat{\lambda}}^j_i}{} \hspace{0.7cm} \downarrow
 \st{{\lambda}^j_i}{}\\
1 \longrightarrow  R_j  \longrightarrow  F_j  \longrightarrow  G_j  \longrightarrow  1.\\
\end{array}\]
The commutativity of this diagram, implies that the homomorphism
${\widehat{\lambda}}^j_i$ maps $R_i$ into $R_j$ and so
${\widehat{\lambda}}^j_i (R_i\cap V(F_i))\subseteq R_j\cap
V(F_j)$. Hence if ${\widehat{\lambda}}^j_i$ maps $S_i$ into $S_j$,
we will have the
 following induced homomorphism, for any $i\leq j$:
\[{\widetilde{\lambda}}^j_i:G_i^*=\frac{F_i}{S_i}\longrightarrow
G_j^*=\frac{F_j}{S_j},\] which forms the directed system $\{G_i^*;
{\widetilde{\lambda}}_i^j,I\}$, called an {\it induced directed
system of covering groups}.

Note that in general, any family of covering groups of a directed
system of groups is not necessarily an induce one. For example we
consider the group ${\bf Z}_2\times {\bf Z}_2$ with two
non-isomorphic covering groups $D_8$ and $Q_8$. So it takes the
trivial directed system $\{G_i; \lambda_i^j,{\bf N}\}$ which
$G_i={\bf Z}_2\times {\bf Z}_2$ and $\lambda_i^j$ to be identity,
for any $i,j\in {\bf N}$. Now if we take the family of covering
groups $\{G_i^*;\ i\in {\bf N}\}$ such that $G_{2i}^*=D_8$ and
$G_{2i+1}^*=Q_8$. Then $\{G_i^*;\ i\in I\}$ does not form induced
directed system.\\
\textbf{\bf Theorem 3.5.}\textit{ Suppose that $\cal V$ is a
Schur-Baer variety. If $\{G_i; \lambda_i^j,I\}$ is a directed
system of finite groups with an induced system of ${\cal
V}$-covering groups $\{G_i^*; {\widetilde{\lambda}}_i^j,I\}$, as
we mentioned above, then the group
$G^*=\displaystyle{\lim_{\longrightarrow}G_i^*}$ is a ${\cal V}$-covering group for
$G=\displaystyle{\lim_{\longrightarrow}G_i}$.}\\
\textit{Proof.} Using the isomorphism $(*)$ and Lemma 3.4, we have
$$ \lim_{\longrightarrow}\frac {R_i}{[R_iV^*F_i]}\cong\lim_{\longrightarrow}\frac
{R_i\cap V(F_i)}{[R_iV^*F_i]}\times \lim_{\longrightarrow}\frac
{S_i}{[R_iV^*F_i]}\ .\ \ \ \ (1)$$ Also by Lemma 2.8(ii), we have
$$ \lim_{\longrightarrow}\frac {R_i}{[R_iV^*F_i]}=\frac
{\displaystyle{\lim_{\longrightarrow}R_i}}
{\displaystyle{\lim_{\longrightarrow}[R_iV^*F_i]}},\ \ \ \ (2)$$
$$\lim_{\longrightarrow}\frac {R_i\cap V(F_i)}{[R_iV^*F_i]}=\frac {\displaystyle{\lim_{\longrightarrow}
(R_i\cap V(F_i))}}{\displaystyle{\lim_{\longrightarrow}[R_iV^*F_i]}}
\ \ \ and\ \ \ \lim_{\longrightarrow}\frac {S_i}{[R_iV^*F_i]}=\frac
{\displaystyle{\lim_{\longrightarrow}S_i}}{\displaystyle{\lim_{\longrightarrow}[R_iV^*F_i]}}\
.\ \ \ \ (3)$$ Therefore using Lemma 2.12 and (1), (2), (3) we
conclude that
\[\displaystyle{\frac
{(\displaystyle{\lim_{\longrightarrow}R_i)}}{[(\displaystyle{\lim_{\longrightarrow}R_i)V^*(\lim_{\longrightarrow}F_i)}]}
\cong \frac {(\displaystyle{\lim_{\longrightarrow}R_i)\cap
V(\lim_{\longrightarrow}F_i)}}{[(\displaystyle{\lim_{\longrightarrow}R_i)V^*(\lim_{\longrightarrow}F_i)}]}
\times \frac
{(\displaystyle{\lim_{\longrightarrow}S_i)}}{[(\displaystyle{\lim_{\longrightarrow}R_i)V^*(\lim_{\longrightarrow}F_i)]}}.}\]
Now by the above relation, Corollary 2.8 and Lemma 3.1,
$\displaystyle{\lim_{\longrightarrow}F_i/\lim_{\longrightarrow}S_i}$
and so the group
\[ \displaystyle{G^*=\lim_{\longrightarrow}G_i^*=\lim_{\longrightarrow}\frac {F_i}{S_i}=
\frac{\displaystyle{\lim_{\longrightarrow}F_i}}{\displaystyle{\lim_{\longrightarrow}S_i}}}\]
will be a $\cal V$-covering group of
$\displaystyle{G=\lim_{\longrightarrow}G_i}$.\ $\Box$

We end this section, by an example showing that the condition of
being directed for index set of the direct system in our study, is essential.\\
\textbf{\bf Example 3.6.} If we omit the condition of being
directed, then the free product of any two groups A, B as a
particular direct limit of groups which its index set is not
directed (see [20]), should have a covering group with the structure
$A^**B^*$, where $A^*$ and $B^*$ are covering groups of A and B,
respectively. But this is a contradiction, when we choose A with
nontrivial Schur multiplier. Since in this case, if $A^**B^*$ is a
covering group of $A*B$, then by Definition 2.4, we will have
\[M(A*B)\cong N\ {\it with}\ N\subseteq Z(A^**B^*)\cap(A^**B^*)'=1.\]
But using a result of Miller [13], we have $M(A*B)\cong M(A)\times
M(B)\neq 1$, which is a contradiction.\\

\begin{center}
\textbf{4. Applications}\\
\end{center}

As we mentioned in the introduction, we have the structure of
covering or generalized covering groups for some famous products
of finitely many groups, such as finite direct and finite regular
products of finite groups [4,24]. In this section, as an
important application of the main result, we present a generalized
covering group for the above products when their index sets are
arbitrary. Also the main result of this note may have an
application in the sense that if one wants to find a generalized
covering group for an arbitrary group, one only needs to find a
generalized covering group for every finitely generated subgroups
of it.

Firstly, by a result of Schur [22], for a finite nilpotent group G
and it's all Sylow subgroups $S_1, S_2,\ldots ,S_n$, we have
$$M(G)\cong M(S_1)\times M(S_2)\times \cdots \times M(S_n).$$
Using this property and also the definition of covering group,
we deduce the following straightforward point.\\
\textbf{\bf Corollary  4.1.}\textit{ Let $G$ be a finite
nilpotent group with it's all Sylow subgroups $S_1, S_2,\ldots
,S_n$, and let $S_i^*$ be a covering group for $S_i$. Then the
group  $G^*= S_1^*\times S_2^*\times \cdots \times S_n^*$ is a
covering group of G.}

Now using the main result of this paper, we conclude the
generalization of the above note, as follows.\\
\textbf{\bf Corollary  4.2.}\textit{ Let $G$ be a torsion
nilpotent group with it's all Sylow subgroups $S_i$, for $i\in
I$. Suppose $S_i^*$ is a covering group of $S_i$, for all
$i\in I$, then the group  $G^*= \prod^{\times}_{i\in I}S_i^*$ is a covering group for G.}\\
\textit{Proof.} First, note that every torsion nilpotent group, is
the direct product of it's Sylow subgroups and so
$G=\prod^{\times}_{i\in I}S_i$. Now if we consider the system
$\{\prod^{\times}_{j\in J}S_j;{\lambda}_J^K\}$, where $J\subseteq
K$ are finite subsets of I and ${\lambda}_J^K$ is the inclusion
map, then the group $G=\prod^{\times}_{i\in I}S_i$ is the direct
limit of this system which is obviously a directed system.

Clearly we have the induced directed system on covering groups of
$\prod^{\times}_{j\in J}S_j$'s with the morphisms
${\widetilde{\lambda}}^K_J$, as follows:
\[\{\prod^{\times}_{j\in J}S_j^*\ ;\ {\widetilde{\lambda}}_J^K\}_{J\st{finite}{\subseteq}I}.\]
Note that the morphisms  ${\widetilde{\lambda}}^K_J$ are clearly
inclusion maps and so the direct product $G^*= \prod^{\times}_{i\in
I}S_i^*$ will be considered as the direct limit of this system and
hence, using Theorem 3.5, the proof is completed.\ $\Box$

Also by a result of G. Ellis [3], for a finite nilpotent group G
and it's all Sylow subgroups $P_1, P_2,\ldots ,P_n$, we have
$${\cal N}_cM(G)\cong {\cal N}_cM(S_1)\times {\cal N}_cM(S_2)\times \cdots
\times {\cal N}_cM(S_n).$$ Using this property and similar
arguments, we deduce the similar corollary for an ${\cal
N}_c$-covering group of a torsion nilpotent
group as follows:\\
\textbf{\bf Corollary  4.3.}\textit{ Let $G$ be a torsion
nilpotent group with it's all Sylow subgroups $P_i$, for $i\in
I$. Suppose $P_i^*$ is an ${\cal N}_c$-covering group of $P_i$,
for all $i\in I$. Then the group  $G^*= \prod^{\times}_{i\in
I}P_i^*$ is an ${\cal N}_c$-covering group for G.}

We next establish the structure of a covering group for the direct
and regular products of arbitrary many of groups which
are generalizations of results of J. Wiegold [24] and W. Haebich [4].\\
\textbf{\bf Corollary  4.4.}\textit{ Let $\{A_i\}_{i\in I}$ be an
arbitrary family of finite groups and suppose that $A_i^*$ is a
covering group of $A_i$, for any $i\in I$. Then the second
nilpotent product of $A_i$'s, is a covering group of
$\prod^{\times}_{i\in I}A_i$.}\\
\textit{Proof.} First, we recall that the second nilpotent product
of a family of  groups as $\{A_i^*\}_{i\in I}$, is defined to be
$$\prod^{\st{2}{*}}_{i\in I}A_i^*=\frac{\prod^{*}_{i\in I}A_i^*}{\gamma_{3}(\prod^{*}_{i\in
I}A_i^*)\cap[A_i^*]^*},$$ where the subgroup $[A_i^*]^*$ is the
kernel of the natural epimorphism \\ $\pi:\prod^{*}_{i\in
I}A_i^*\rightarrow \prod^{\times}_{i\in I}A_i^*$, which has also
the following structure:
$$[A_i^*]^*=\langle\ [A_i,A_j]\ ; i\neq j\ \rangle^{\prod^{*}_{i\in I}A_i^*}.$$
Now similar to the proof of Corollary 4.2, we consider the directed
system $\{\prod^{\times}_{j\in J}A_j\ ;\ {\lambda}_J^K\}$ of finite
direct products and the inclusion maps, with
 the group $A=\prod^{\times}_{i\in I}A_i$ as the direct limit of this system.

By a result of J. Wiegold [24], the group $\prod^{\st{2}{*}}_{i\in
J}A_i^*$ is a covering group of $\prod^{\times}_{j\in J}A_j$,
where $J$ is finite. Consider $\{\prod^{\st{2}{*}}_{i\in J}A_i^*\
;\ {\widetilde{\lambda}}_J^K\}$ as a directed system, where the
morphism ${\widetilde{\lambda}}^K_J$ is clearly the inclusion map,
${\widetilde{\lambda}}_J^K:\prod^{\st{2}{*}}_{j\in
J}A_j^*\hookrightarrow \prod^{\st{2}{*}}_{k\in K}A_k^*$. Clearly
$A^*=\prod^{\st{2}{*}}_{i\in I}A_i^*$ is the direct limit of the
last induced directed system and so using Theorem 3.5 is a
covering
group for $\prod^{\times}_{i\in I}A_i$.\ $\Box$\\
\textbf{\bf Notation and Corollary  4.5.}\textit{ Let $G$ be a
regular product of it's finite subgroups $A_i$, $i\in I$, where
$I$ considered as an ordered set. For each $i\in I$, $L_i$
denotes a fixed covering group for $A_i$ and consider the exact
sequence $1\ra M_i\ra L_i \st {\nu_i}{\ra} A_i\ra 1$ such that
$$ M_i\subseteq Z(L_i)\cap L'_i\ \ and\ \ M_i\cong M(A_i),$$}
\textit{where $M_i$ is a normal subgroup of $L_i$, and $M(A_i)$ is
the Schur multiplier of $A_i$.} \textit{Assume that
$A=\prod_{i\in I}^{*}A_i$ and $L=\prod_{i\in I}^{*}L_i$ are free
products of $A_i$'s and $L_i$'s, respectively. We denote by $\nu$
the natural homomorphism from $L$ onto $A$ induced by the
$\nu_i$'s. Also, if $\psi$ is the natural homomorphism from $A$
onto $G$ induced by the identity map on each $A_i$,
$$ L=\prod_{i\in I}^*L_i\st {\nu }{\lra }A=\prod_{i\in I}^*A_i\st {\psi}{\lra}
G\ra 1\ ,$$ then we denote by $H$ the kernel of $\psi$ and set}

\[J=\nu^{-1}(H)\cap[L_i^L],\ \ \ N=\prod_{\st{i,j=1}{i\neq j}}^{n}[M_i,L_j]^L,\ \ \
M=(\prod_{i=1}^{n}M_i)J.\]
\textit{ Finally, $\bar{L}$ and
$\bar{M}$ denote the images of $L$ and $M$ under the natural
homomorphism $L\rightarrow L/N[J,L]$, respectively. Then there is
an exact sequence}
$$1 \rightarrow \bar{M} \rightarrow \bar{L}
\rightarrow G \rightarrow 1, $$
\textit{such that $\bar{M}
\subseteq Z(\bar{L})\cap [\bar{L},\bar{L}]$ and $\bar{M}\cong
M(G)$. In particular, $\bar{L}$ is a covering group
of $G$.}\\
\textit{Proof.} First of all, we note that the group $G$ is called
the regular product of it's subgroups $A_i$'s, with the ordered
set $I$, if the following two conditions hold:
$$G=\langle A_i\ ;\ i\in I\rangle\ \ \ \ ,\ \ \ \ A_i\cap \widehat{A}_i=1\ (\forall i\in I),$$
where $\widehat{A}_i$ is defined to be the group
$$\widehat{A}_i=\prod_{j\in J,j\neq i}A_j^{G}.$$

Now similar to the previous notes and using the definition, we
clearly consider the directed system on finite regular products of
$A_i$'s, with inclusion maps. As we saw before, it induces a
directed system on their covering groups, with induced
homomorphisms which are also inclusion.

By a result of W. Haebich [4], the corollary holds for any finite
index set $I$. Now, using Theorem 3.5, it is easy to check that
the group $\bar{L}$ as a direct limit of the induced system is a
covering group of the regular products of $A_i$'s which is
considered as a
direct limit of the first system.\ $\Box$\\

Note that the above theorem is in fact a generalization of the
Haebich's formula for a covering group of any finite regular product
of finite groups [4]. However the main proof of Haebich could be
easily generalized to the regular product of infinitely many of
finite groups, but our proof as an application of the main result is
another proof to this generalization.\\
{\large \bf Acknowledgment}

The authors would like to thank the referee for giving attention to
the paper and spending a good amount of time.

This research was in part supported by a grant from Center of
Excellence in Analysis on Algebraic Structures, Ferdowsi University
of Mashhad.

\end{document}